\documentclass[12pt,bezier]{article}
\usepackage{amsmath,amssymb,amsfonts}

\newcommand{\D}{\displaystyle}
\newcommand{\T}{\mathrm}
\newcommand{\lo}{\longrightarrow}
\newcommand{\ov}{\overline}
\newcommand{\sub}{\subseteq}
\newcommand{\st}{\stackrel}

\begin{document}
\begin{center}
{\bf\large Ljusternik-Schnirelmann Categories,}\\ {\bf\large Links
and Relations }\\ [1cm] M.R. Razvan \\  \noindent {\it \small
Department of Mathematics, Sharif University of Technology\\
P.O.Box: $11365-9415$, Tehran, IRAN
\\ Institute for Studies in Theoretical
Physics and Mathematics\\  P.O.Box: $19395-5746$,
Tehran, IRAN}\\ E-Mail:
 razvan@karun.ipm.ac.ir \\ Fax: 009821-2290648

\end{center}

\vspace*{5mm}
\begin{abstract}
This paper is concerned with some well-known
Ljusternik-Schnirelmann categories. We desire to find some links
and relations among them. This has been done by using the concepts
of precategoty, T-collection and closure of a category.
\end{abstract}


\section{Introduction}

In Ljusternik-Schnirelmann theory, there are several categories
with different properties. First of all, Ljusternik introduced a
category by means of closed contractible sets which is known as
(classical) Ljusternik-Schnirelmann category \cite{L}. Then
Schnirelmann estimated this category by using cohomology theory.
This method is still the most convenient way to compute LS
categories. For this reason, some authors use the same notations
for cohomological category \cite{MS}. The cuplength category is
also derived from this approach to LS theory. Besides these
categories, there is another LS category which is similar to the
classical one and it is defined by means of open contractible
subsets \cite{J}. For the excision property of open subsets, this
category is preferable especially in computation.

In the literature, all of the above categories are known as
Ljusternik-Schnirelmann category and it is customary to use the
notation $cat$ for them. In order to distinguish these categories,
we call them: (classical) Ljusternik-Schnirelmann (shortly LS)
category, Cohomology Ljusternik-Schnirelmann (CLS) category, Cup
Length (CL) category and Homotopy Ljusternik-Schnirelmann (HLS)
category, respectively. The aim of this paper is to find some
links and relations among these categories.

First of all in section 2, we define the concepts of category and
precategory. Indeed, both LS and CL categories are precategories
(but not necessarily categories) in an arbitrary topological
space. In section 3, we compare LS and HLS categories and find a
relation between them. In section 4, we introduce the concept of
T-collection which gives an extension of HLS and CLS categories.
Finally in section 5, we prove that the cuplength map is a
precategory and hence defines a category. This category has been
used in \cite{CZ} to prove the Arnold conjecture for even
dimensional tori. We also find a relation between CL and CLS
categories in this section.

\section{Category}
Let $M$ be a topological space. A category on $M$ is a map $\nu:
2^M\lo \mathbb{Z} \cup \{+\infty\}$ which satisfies the following
axioms:

i) If $A\subset B$, then $\nu(A)\leq \nu(B)$.

ii) $\nu (A\cup B)\leq \nu (A)+\nu (B)$.

iii) For every subset $A\subset M$ there exists an open set $U\subset M$ such
that $A\subset U$ and $\nu(A)=\nu(U)$.

iv) If $f:M\lo M$ is a continuous map homotopic to the identity $id_M$, then
$\nu(A)\leq \nu (f(A))$ for every subset $A\subset M$.\\

It has been shown that if $\varphi ^t$ is a continuous flow on a
compact metric space $M$ with a Morse decomposition {$M_1,\cdots
,M_n$}, then $\nu (M) \leq \sum _{i=1} ^n \nu(M_i)$ for every
category $\nu$ on $M$. This result leads to a critical point
theory for gradient-like flows. (See \cite{CZ} and \cite{R} for
details.)

\noindent {\bf Example 1.} (Ljusternik-Schnirelmann category). A
metric space $X$ is an absolute neighborhood extensor, shortly an
ANE, if for every metric space $Y$, every closed subset $B$ of $Y$
and $f:B\lo X$ continuous, there exists a continuous extension of
$f$ defined on a neighborhood of $B$
 in $Y$. For a subset $A\subset X$, we define $\nu_{LS}(A)$ to be the
minimum number of closed sets contractible in $M$ required to cover $A$. When
$X$ is an ANE, $\nu_{LS}$ is a category on $X$ (See \cite{MW}, $\S$ 4.6.)
which satisfies the following axiom:

v) If A consists of a single point, then $\nu (A)=1$.

Note that $\nu_{LS}$ is defined on  every topological space and satisfies
axioms
(i), (ii), (iv) and (v). But there are examples of compact metric spaces in
which $\nu_{LS}$ does not satisfy axiom (iii). \\

\noindent {\bf Example 2.} (Precategory) Let $X$ be a topological
space and $\tau$ be the topology on $X$. A precategory is a map
$\nu_0: \tau\lo \mathbb{Z} \cup\{+\infty\}$ which satisfies axioms
(ii) and (iv). For example $\nu_{LS}$ defines a precategory on
every topological space. In section 5, we shall show that the
cuplength map is also a precategory. Given a precategory $\nu_0$
on $X$, we can define a category $\st{\sim}{\nu_0}: 2^M \lo
\mathbb{Z} \cup \{+\infty\}$ as follows: $$\st{\sim}{\nu_0}
(A)=\min \{\nu_0 (U) | U \subset X \ \T{open, containing}
 \ A\}$$
Notice that $\nu_{LS}$ defines a category on $X$ if and only if
$\nu_{LS}=\st{\sim}{\nu}_{LS}$. Thus $\st{\sim}{\nu}_{LS}$ can be
considered as a generalization of Ljusternik-Schnirelmann
category.
\\

\noindent {\bf Example 3.} (Closure of a category) Let $\nu$ be a
category on a normal topological space $X$. We define
$\ov{\nu}:2^M\lo \mathbb{Z} \cup\{\infty\}$ by
$\ov{\nu}(A)=\nu(\ov{A})$. It is not hard to see that $\ov{\nu}$
is a category. (Normality is used to prove axiom (iii)). We have
$\nu\leq \ov{\nu}$ and if $\nu$ satisfies axiom (v), then so does
$\ov{\nu}$.

\section{ HLS Category}

Let $M$ be a topological space. A subset $A\subset M$ is called
contractible in $M$ if the inclusion map $A\lo M$ is homotopic to
a constant. We define Homotopy Ljusternik-Schnirelmann category as follows:\\

\noindent {\bf Definition.}
 The HLS-category $\nu_H(A)=\nu_H(A,M)$ of a subset $A\subset M$ is
defined to be the minimum number of open sets contractible in $M$
required to cover $M$. If such a covering does not exist, we set
$\nu_H(A)=+\infty$ and if it exists, $A$ is called
$H$-categorizable (in $M$).

It is easy to see that $\nu_H$ satisfies axioms (i)-(iii).
Moreover a subset $A\subset M$ is H-categorizable if and only if
$\nu_H(\{x\})=1$  for every $x\in A$. Thus
 $\nu_H$ satisfies axiom (v) if and only if $M$ is $H$-categorizable.
 The following lemma gives a generalization of axiom (iv). \\

\noindent {\bf Lemma 3.1.} If $Y$ dominates $X$, i.e. there are
continuous maps
 $f:X\lo Y$ and
 $g:Y\lo X$ with  $g\circ f\sim
id_X$, then for every $H$-categorizable subset $A\subset Y$,
$f^{-1} (A)$ is $H$-categorizable and $\nu_H(f^{-1}(A))\leq
 \nu_H(A)$. In particular if $Y$ is $H$-categorizable, then
so is $X$ and $\nu_H(X)\leq \nu_H(Y)$. Hence $\nu _H$ is an
invariant of homotopy type. \\ {\bf Proof.} It is enough to prove
that for every open set $U\subset Y$ contractible in $Y$,
$f^{-1}(U)$ is contractible in $X$. Consider the following
commutative diagram in which $i$ and $j$ are inclusion maps:

\unitlength=1mm \special{em:linewidth 0.4pt} \linethickness{0.4pt}
\begin{picture}(95.33,43.33)
\put(80.00,40.00){\vector(-1,0){20.00}}
\put(50.00,35.00){\vector(0,-1){10.00}}
\put(90.00,35.00){\vector(0,-1){10.00}}
\put(80.00,20.00){\vector(-1,0){20.00}}
\put(70.33,23.67){\makebox(0,0)[cc]{$(f|_{A})^*$}}
\put(70.33,43.33){\makebox(0,0)[cc]{$f^*$}}
\put(50.00,40.00){\makebox(0,0)[cc]{$H^k(X)$}}
\put(46.67,20.00){\makebox(0,0)[cc]{$H^k(A)$}}
\put(54.00,30.00){\makebox(0,0)[cc]{$i^*$}}
\put(95.33,30.33){\makebox(0,0)[cc]{$j^*$}}
\put(90.00,20.00){\makebox(0,0)[cc]{$H^k(f(A))$}}
\put(90.00,40.00){\makebox(0,0)[cc]{$H^k(Y)$}}
\end{picture}

\vspace*{-1cm}

$f\circ i=j\circ f|_V \Rightarrow g \circ f \circ i = g\circ j
\circ f|_V \Rightarrow i \sim g\circ j \circ f|_V \sim$ constant.
$\square$\\

The rest of this section concerns a comparison between $\nu _H$
and $\nu _{LS}$. We first prove that every ANE is H-categorizable.

\noindent {\bf Definition.}
 A topological space $X$ is called semi-locally contractible if for every
$x\in X$ and open set $U\subset X$ with $x\in U$, there exists a
neighborhood $V$ of $x$ such that $x\in V\subset U$ and $V$ is
contractible in $U$. \\

\noindent {\bf Proposition 3.2.} Every ANE is semi-locally
contractible, hence it is $H$-categorizable. \\ {\bf Proof.}
Suppose that $X$ is an ANE, $x_{0}
 \in X$ and $U\subset X$ is an
open set such that $x_0\in U$. We define $f_0:\ov{U} \times \{0,1\} \cup \{x_0\}
\times [0,1]\lo X$ by $f_0(x,t)=\left\{ \begin{array}{ll}
x & \text{if} \ t=0\\[-0.2cm]
x_0 & \text{otherwise.}
\end{array}\right .$ \\
$f_0$ is continuous and $\ov{U}\times \{0,1\}\cup
\{x_0\}\times[0,1]$ is a closed subset of $X\times [0,1]$. Since
$X$ is an ANE, there exists a continuous extension $f_1$ of $f_0$
defined on a neighborhood of $\{x_0\}\times [0,1]$. We have
$f_1(x_0,t)=f_0(x_0,t)=x_0\in U$. Thus $f_1^{-1}(U)$ is an open
set containing $\{x_0\}\times [0,1]$. By compactness of $[0,1]$,
there is an open set $V\subset U$ such that $V\times [0,1] \subset
f_1^{-1} (U)$. Now consider $f=f_1|_{V\times [0,1]}: V\times
[0,1]\lo U$. We have: $$f(x,0)=f_1(x,0)=f_0(x,0)=x \ \ \text{and}
\ \ f(x,1)=f_1(x,1)=f_0(x,1)=x_0.$$ It means that $V$ is
contractible in $U$. $\square$ \\

\noindent {\bf Proposition 3.3.} If $\nu_{LS}$ defines a category
on a normal topological space $X$, then $\nu_{LS}=\ov{\nu}_H$.
 In particular $\nu_{LS}$ and $\nu_H$ agree on closed
subsets of an $ANE$. \\ {\bf Proof.} Step 1. $\nu_{LS}
(A)=\nu_{LS}(\ov{A})$ for every subset $A\subset X$. \\ Since
$\nu_{LS} (A) \leq \nu_{LS} (\ov{A})$, we may assume that
$\nu_{LS}(A)<\infty$. Let $A\subset \D{\bigcup_{i=1}^{n}} A_i$
where each $A_i$ is closed contractible in $X$. Thus $\ov{A}
\subset \D{\bigcup_{i=1}^n} A_i$ and it follows that $\nu_{LS}
(\ov{A})\leq \nu_{LS}(A)$.

{Step 2.} $\nu_H(A)\leq \nu_{LS}(A)$ for every subset $A\subset
X$.\\ Suppose that $A\subset \D{\bigcup_{j\in J}} A_j$ and each
$A_j$ is closed contractible in $X$. Thus $\nu_{LS}(A_j)=1$ for
every $j\in J$ and by axiom (iii), there exists an open set $U_j$
with $A_j\subset U_j$ and $\nu_{LS}(U_j)=1$. Now each $U_j$ is
open contractible in $X$ and $A\subset \D{\bigcup_{j\in J}} U_j$.

{Step 3.} For every closed subset $A\subset X$, $\nu_{LS}(A)\leq
\nu_H(A)$.\\ We may assume that $\nu_H(A)<\infty$. Suppose that
$A\subset\D{\bigcup_{i=1}^{n}} U_i$ where each $U_i$ is open
contractible in $X$. We have $A=\D{\bigcup_{i=1}^n} A\cap U_i$ and
each $A\cap U_i$ is open in $A$. By Shrinking Lemma \cite{M},
 there are $V_i\subset A$ open in $A$ such
that $A=\D{\bigcup_{i=1}^n}
 V_i$ and $\ov{V}_i \subset A\cap U_i$. (The closure
is taken in $A$ which is a closed subset of $X$ and there is no
ambiguity.) Since $\ov{V}_i\subset U_i$ and $U_i$ is contractible
in $X$, each $\ov{V}_i$ is closed contractible in $X$ and
$A\subset \D{\bigcup_{i=1}^n} \ov{V}_i$.

The proof of Proposition 3.3. is easy now. By steps 1 and 2,
$\nu_H(\ov{A})\leq \nu_{LS} (\ov{A})=\nu_{LS} (A)$ and by steps 1
and 3, $\nu_{LS} (A)=\nu_{LS} (\ov {A})\leq \nu_H(\ov{A})$.
$\square$

the relation $\nu _{LS}=\ov{\nu}_H$ says that $\nu _{LS}$ can be
obtained by $\nu _H$. Moreover the above two results show that
$\nu _H$ is applicable to more spaces than $\nu _{LS}$.

\section{T-collection}
Having a glance at the definition of HLS category, one can easily
find a general method to define a category. We do this by using
the concept of T-collection. A generalization of this concept has
been used in \cite{R} to obtain a result in critical point theory.
\\ {\bf Definition.}
(a) A collection $T$  of open subsets of $M$ is called a
T-collection if for every continuous map $f:M\lo M$ homotopic to
$id_M$, $U\in T$ implies $f^{-1} (U) \in T$. Every $U\in T$ is
called $T$-trivial in $M$. A subset $A\subset M$ is called
$T$-categorizable if $A$ is covered by some $T$-trivial subsets.

(b) For every T-collection $T$ on $M$, we define the associated
category $\nu_T(A)$ to be the minimum number of open sets in $T$
required to cover $A\subset M$. If such a covering does not exist,
we set $\nu_T(A)=+\infty$. It is not hard to see that $\nu_T$ is a
category on $M$. \\

\noindent {\bf Example 1.} Let $U\subset M$ be an open set and $T$
be the set of all $f^{-1}(U)$ where $f:M\lo M$ is continuous
homotopic to $id_M$. Then $T$ is a T-collection.\\ {\bf Example
2.} If $T_j$'s are T-collections for $j\in J$, then so are
$\D{\bigcap_{j\in J}} T_j$ and $\D{\bigcup_{j\in J}} T_j$. It is
easy to see that every T-collection can be obtained by these two
examples. \\

\noindent {\bf Definition.} Suppose that $\nu:2^M\lo \mathbb{Z}
\cup \{+\infty\}$ satisfies axiom (iv). For every $n\in
\mathbb{Z}$ we define $T_{\nu,n}=\{U \subset M | U$ is open and
$\nu(U)\leq n\}$ which is obviously a T-collection. If $\nu$ is a
category, then $T_{\nu,1}$ is called the associated precategory
and it is denoted by $T_{\nu}$. \\

Now there are two natural questions: 1. When a given category is
associated to some T-collection? 2. When a given T-collection is
associated to some category? The rest of this section concerns
these two questions. First of all, notice that by the definition
of $\nu_T$, we have $T\subseteq T_{\nu_T}$. \\

\noindent {\bf Lemma 4.1.} For every category $\nu$ on $M$ and
$n\in\mathbb{N}$, $\nu\leq n\nu_{T_{v,n}}$. In particular $\nu\leq
\nu_{T_{\nu}}$. \\ {\bf Proof.} Let $A$ be any subset of $M$. We
may assume that $\nu_{T\nu}(A)<\infty$ and $A\subset
\D{\bigcup_{i=1}^{m}} U_i$ where $U_i\in T_{\nu, n}$ and
$m=\nu_{T_{\nu,n}}(A)$. Now by axiom (ii), $\nu(A)\leq
\D{\sum_{i=1}^{m}}\nu(U_i)\leq mn=n\nu_{T_{\nu,n}}(A)$.
$\square$\\

\noindent {\bf Proposition 4.2.} i) $T_{\nu_{T_{\nu}}}=T_{\nu}$
for every category $\nu$ on $M$. \\ ii) $\nu_{T_{\nu_T}}=\nu_T$
for every T-collection $T$ on $M$.\\ {\bf Proof.} Notice that if
$\nu \leq \nu'$ then $T_{\nu}\supseteq T_{\nu'}$ and if
$T\subseteq T'$ then $\nu_{T}\geq \nu_{T'}$. Now we have \\
$T\subseteq T_{\nu_T}$, hence $T_{\nu}\subseteq T_{\nu_{T_{\nu}}}
\text{and} \ \nu_{T}\geq \nu_{T_{\nu_T}}$. Similarly $\nu\leq
\nu_{T_\nu}$, hence $\nu_{T_{\nu_T}} \text{and} \ T_{\nu}\supseteq
T_{\nu_{T_{\nu}}}$. $\square$

\noindent {\bf Corollary 4.3} i) For every category $\nu$ on $M$,
$\nu=\nu_{T_{\nu}}$ if and only if $\nu=\nu_T$ for some
T-collection $T$ on $M$. \\ ii) For every T-collection $T$ on $M$,
$T=T_{\nu_{T}}$ if and only if $T=T_{\nu}$ for some category $\nu$
on $M$. \\ {\bf Proof.} If $\nu=\nu_T$, then
$\nu_{T_{\nu}}=\nu_{T_{\nu_{T}}}=\nu_T$ by the above proposition
and hence $\nu_{T_{\nu}}=\nu$. Similarly, if $T=T_{\nu}$, then
$T_{\nu_T}=T_{\nu_{T_{\nu}}}=T_{\nu}=T$. $\square$

\section{Cohomological Category}

This section concerns two examples: Cohomology
Ljusternik-Schnirelmann category and Cup Length category. From now
on, $H^*$ is assumed to be a fixed cohomology functor on a
topological space $M$ with a cup product $\cup$. \\

\noindent {\bf Definition.} A subset $A\subset M$ is called
cohomologically trivial in $M$ if the restriction map $j^*:
H^k(M)\lo H^k(A)$ is zero for $k>0$. The set of all
cohomologically
 trivial open subsets of $M$ is denoted by $T_c(M)$. For every
subset $A\subset M$, we define the cup length of $A$ (in $M$) to
be the minimum integer $N>0$ such that for any set of cohomology
classes $\alpha_j\in H^{k_j} (M)$, $j=1,\cdots , N$ of degree
$k_j\geq 1$, the class $(\alpha_1 \cup \cdots \cup
\alpha_N)|_A=j^* (\alpha_1 \cup \cdots \cup \alpha_N)$ is zero.
Hence $A$ is cohomologically trivial if and only if
$cuplength(A)=1$. \\

\noindent {\bf Proposition 5.1.} Suppose that $A$ and $B$ are
excisive in $X$ and $\alpha,\beta \in H^* (X)$ with
$\alpha|_A=\beta|_B=0$, then $(\alpha \cup \beta)|_{A\cup B}=0$.
In particular if $X=U \cup V$ where $U$ and $V$ are open sets and
$\alpha|_U=\beta|_V=0$, then $\alpha \cup \beta=0$.\\ {\bf Proof.}
Since $A$ and $B$ are excisive, the cup product map
$H^*(X,A)\otimes H^*(X,B)\lo H^* (X,A\cup B)$ is defined \cite{D}
and the following diagram commutes: $$\begin{array}{ccccc}
H^*(X,A) & \otimes & H^*(X,B) & \lo & H^* (X, A\cup B)\\
\downarrow & & \downarrow & & \downarrow\\ H^*(X) & \otimes &
H^*(X) & \lo & H^*(X)\\ \downarrow & & \downarrow & & \downarrow\\
H^*(A) & \otimes & H^*(B) &  & H^*(A\cup B)
\end{array}$$
the rest of the proof is a diagram chase and similar to Lemma 4 in
\cite{CZ}. $\square$

\noindent {\bf Lemma 5.7.} Suppose that $f:X\lo Y$ is a continuous
map such $f^*:H^k(Y)\lo H^k(X)$ is surjective for $k>0$. Then
$cuplength(A)\leq$ $cuplength(f(A))$ for every $A\sub X$, . \\
{\bf Proof.} We may assume that $m=cuplength(f(A))<\infty$. Now
for $\alpha_i \in H^{k_i} (X), i=1,\cdots , m$ with $k_i>0$, there
are $\beta_i\in H^{k_i} (Y)$ such that $f^*(\beta_i)=\alpha_i$.
Consider the following diagram

\unitlength=1mm \special{em:linewidth 0.4pt} \linethickness{0.4pt}
\begin{picture}(95.33,43.33)
\put(80.00,40.00){\vector(-1,0){20.00}}
\put(50.00,35.00){\vector(0,-1){10.00}}
\put(90.00,35.00){\vector(0,-1){10.00}}
\put(80.00,20.00){\vector(-1,0){20.00}}
\put(70.33,23.67){\makebox(0,0)[cc]{$(f|_{A})^*$}}
\put(70.33,43.33){\makebox(0,0)[cc]{$f^*$}}
\put(50.00,40.00){\makebox(0,0)[cc]{$H^k(X)$}}
\put(46.67,20.00){\makebox(0,0)[cc]{$H^k(A)$}}
\put(54.00,30.00){\makebox(0,0)[cc]{$i^*$}}
\put(95.33,30.33){\makebox(0,0)[cc]{$j^*$}}
\put(90.00,20.00){\makebox(0,0)[cc]{$H^k(f(A))$}}
\put(90.00,40.00){\makebox(0,0)[cc]{$H^k(Y)$}}
\end{picture}

\vspace*{-2cm}

\noindent where $i^*$ and $j^*$ are restriction maps. Now we have
$(\alpha_1 \cup \cdots \cup \alpha_m)|_A=i^*(\alpha_1 \cup\cdots
\cup \alpha_m)=i^*(f^*(\beta_1)\cup\cdots \cup (f^*(\beta_m))=
i^*(f^*(\beta_1\cup \cdots \cup \beta_m))=0$ and the proof in
complete. $\square$

The above two propositions show that the cup length
 map is a precategory and hence it
defines a category which is denoted by $\nu_{CL}$. Moreover if we
use the above lemma in the case of $cuplength(A)=1$, we see that
$T_c$ is a T-collection and hence defines a category which is
denoted by $\nu_c$. It is easy to see that
$T_{\nu_{CL}}=T_c\cup\{\varnothing\}$, hence $\nu _
c=\nu_{T_{\nu_{CL}}}$ which is a relation between CLS and CL
categories. Moreover by Lemma 4.1, $\nu_{CL}\leq
\nu_{T_{\nu_{CL}}}=\nu_c$ and since $T_H\subset T_c$, we have
$\nu_{CL} \leq \nu_c \leq \nu_H$. These inequalities show that the
best estimate is provided by $\nu_H$, but $\nu_{CL}$ is somewhat
easier to compute.

\paragraph{Acknowledgments.} The author  wishes to thank professors S. Shahshahani
and M. Hesaaraki for helpful discussions.


\begin{thebibliography}{1}

\bibitem{CZ} Conley, C. and Zehnder, E., The Birkhoff-Lewis fixed
point theorem and a conjecture of V.I. Arnold. {\it Inventions
 Mathematicae}, 73 (1983), 33-49.


\bibitem{D} Dold, A., Lectures on Algebraic Topology,
 Springer-Verlag, Berline, Heidelberg, 1972.

\bibitem{J} James, I.M., On category in thge sense of
Lusternik-Schnirelman, {\it Topology} 17 (1978), 331-248.

\bibitem{L} Ljusternik,L., Topology and the calculus of
variations, {\it Uspehi. Mat. Nauk} 1(11) no.1, (1946), 30-56.
MR9, 51.

\bibitem{M} Munkres, J.R., Topology, A First Course, Prentice-Hall Inc.
Englewood Cliffs, New Jersy, 1975.

\bibitem{MS} McDuff, D. and Salamon, D., Introduction to
Symplectic Topology, Clarendon Press, Oxford, 1995.

\bibitem{MW} Mawhin, J. and Willem, M., Critical Point Theory and
Hamitonian Systems, Springer-Verlag, New York, Inc., 1989.

\bibitem{R} Razvan, M.R., Ljusternik-Schnirelmann cateories on index pairs,
 to appear in {\it Math. Z.} 
\end{thebibliography}
\end{document}